\documentclass[letterpaper,10 pt,conference]{ieeeconf}
\IEEEoverridecommandlockouts  

\usepackage{amsmath,amssymb,euscript ,yfonts,psfrag,latexsym,dsfont,graphicx,bbm,color,amstext,wasysym,subfig,parskip}

\graphicspath{{./},{./figures/},{../figures/}}

\newtheorem{thm}{Theorem}
\newtheorem{theorem}{Theorem}
\newtheorem{cor}[thm]{Corollary}

\newcommand{\cG}{\mathcal G}
\newcommand{\cL}{\mathcal L}
\newcommand{\cV}{\mathcal V}

\newcommand{\diag}{\operatorname{diag}}
\newcommand{\RB}{{\,^{\rm RB}}}
\newcommand{\one}{{\mathbbm 1}}
\newcommand{\zz}{{0\hspace*{-6pt}0}}

\newcommand{\cE}{\mathcal{E}}

\newcommand{\cS}{\mathcal{S}}

\newcommand{\bes}{\begin{equation*}}
\newcommand{\ees}{\end{equation*}}
\newcommand{\beas}{\begin{eqnarray*}}
\newcommand{\eeas}{\end{eqnarray*}}
\newcommand{\bea}{\begin{eqnarray}}
\newcommand{\eea}{\end{eqnarray}}
\newcommand{\be}{\begin{equation}}
\newcommand{\ee}{\end{equation}}

\newcommand{\bbl}{\begin{block}}
\newcommand{\ebl}{\end{block}}

\definecolor{grey}{rgb}{0.6,0.6,0.6}
\definecolor{lightgray}{rgb}{0.97,.99,0.99}

\begin{document}
\title{Ruelle-Bowen continuous-time random walk}
\author{Yongxin Chen, Tryphon T. Georgiou, Michele Pavon
\thanks{Y.\ Chen Department of Electrical and Computer Engineering,
  Iowa State University,
  Ames, IA 50011; email: yongchen@iastate.edu}
\thanks{T.\ T. Georgiou is with the Department of Mechanical and Aerospace Engineering, University of California, Irvine, CA; email: tryphon@uci.edu}
\thanks{M. Pavon is with the Dipartimento di Matematica ``Tullio Levi Civita",
Universit\`a di Padova, 35121 Padova, Italy; email: pavon@math.unipd.it}}

\maketitle

\begin{abstract}
We define the probability structure of a continuous-time time-homogeneous Markov jump process, on a finite graph, that represents the continuous-time counterpart of the so-called Ruelle-Bowen discrete-time random walk. It constitutes the unique jump process having maximal entropy rate. Moreover, it has the property that, given the number of jumps between any two specified end-points on the graph, the probability of traversing any one of the alternative paths that are consistent with the specified number of jumps and end-points, is the same for all, and thereby depends only on the number of jumps and the end-points and not the particular path being traversed.
\end{abstract}

\section{Introduction}
The motivation for this note stems from our recent work \cite{chen2017robust,chen2017thermodynamics}, where the topic was the scheduling for transporting resources over a graph. Indeed, a theoretical framework was developed based on the Schr\"odinger bridge problem, namely, to identify a probability law on alternative paths through which (probability) mass is transported between end-point marginals. The sought probability law was taken as the closest to that of the discrete-time Ruelle-Bowen (RB) random walk \cite{DelLib11,parry1964intrinsic,Rue04}, so as to disperse the flow of mass maximally over available paths and, thereby, ensure a level of robustness.

In the present note, we discuss and develop the continuous-time counterpart of the discrete-time RB-walk, as a preparation towards subsequent development of transport protocols emulating our previous aforementioned work \cite{chen2017robust,chen2017thermodynamics}. The construction builds on a standard model for jump processes as well as the discrete-time RB-probability law.

\section{Markov jump processes}\label{sec:realization}

We consider a time-homegenious Markov jump process $\{X_t\mid t\geq 0\}$ with finite state-space $\cV=\{1,2,\ldots,n\}$. This is identified with the vertices (nodes) of a strongly connected directed graph $\cG=(\cV,\cE)$ having no self loops. We suppose that the transition probabilities $p_{ij}(t)$ are continuous at $t=0$ and that $\lim_{t\to 0}p_{i,j}(t) =\delta_{ij}$, where $\delta_{ij}=1$ when $i=j$ and $0$ otherwise. Then, as is well known, the transition probabilities are of the form
\[
\left[ p_{ij}(t)\right]_{i,j=1}^n= \exp(Qt),
\]
where the infinitesimal generator 
\[Q=\left[ q_{ij}\right]_{i,j=1}^n\]
of the process satisfies
\begin{align*}
-q_{ii}&=:q_i>0,\\
q_{ij}&\geq 0 \mbox{ for }i\neq j,\\
q_i&=\sum_{j\neq i}q_{ij}.
\end{align*}
That is, $Q\one=\zz$, where $\one, \zz$ denote column vectors with ones and zeros, respectively. We note that the chain is ergodic, $p_{ij}(t)>0$ for all $i,j$ and $t>0$, which follows from the fact that the graph is strongly connected.

It is quite standard to realize the jump process \cite[Chapter 2]{norris} via the bivariate process
\[
\{(Z_k,T_k) \mid k=0,1,\ldots\}
\]
on the product space $\cS=\cV\times (0,\infty)$, where $Z_k$ is a discrete-time Markov chain with transition probabilities
\[
\pi_{ij} ={\mathbb P}\{Z_{k+1}=j\mid Z_{k}=i\} = \frac{q_{ij}}{q_i},
\]
taking place at a random time $t_{k+1}=\sum_{i=0}^kT_i, t_0 = 0$, where $T_k$ are random times exponentially distributed with density $q_{\,_{\hspace*{-1pt}Z_k}}\exp(-q_{\,_{\hspace*{-1pt}Z_k}}t)$ representing time-intervals between successive transitions between elements of $\cV$.
Thence, $X_t$ is the right-continuous process
\[
X_t = Z_k \mbox{ for } t_k\leq t<t_{k+1}.
\]
Note that, in general, the random time $T_{k}$ is a function of $Z_k$. The bivariate realization of $X_t$ models a random walker on the network, taking steps and transitioning from one node to the next, according to this probabilistic model.

\section{Discrete-time Markov chains and the Ruelle Bowen random walk}\label{sec:discrete}

We briefly outline the construction of the discrete-time RB random walk \cite{DelLib11} on $\cG$, as it will be needed in the continuous-time construction that follows. 

Once again we begin with a strongly connected directed graph\footnote{In discrete-time, self loops are allowed.} $\cG=(\cV,\cE)$, having nodes $\cV=\{1,2,\ldots,n\}$, whose adjacency matrix is denoted by $A=[a_{ij}]_{i,j=1}^n$. We consider a stationary, i.e., time-homogeneous, discrete-time Markov chain $\{Z_k \mid k=0,1,\ldots\}$ taking values on $\cV$, and we let $P=[p_{ij}]_{i,j=1}^n$ denote the matrix of transition probabilities and $\pi=\left[\pi_i\right]_{i=1}^n$ the corresponding stationary probability distribution.\footnote{Vectors, such as $\pi$, will be thought as column vectors throughout.} The entropy rate of the process is
    \begin{equation}\label{eq:entropyd}
        H(P)=-\sum_{i,j} \pi_i p_{ij} \log p_{ij}.
    \end{equation}
The discrete-time RB random walk is the Markov chain that is consistent with the topology of the graph and has maximal entropy rate. It is well known that this is unique \cite{parry1964intrinsic,DelLib11} and the construction proceeds as follows. Its transition probability matrix $P^\RB=[p_{ij}]_{i,j=1}^n$ has entries
\[
p^\RB_{ij}=\frac{1}{\lambda_A}\frac{\varphi_j}{\varphi_i},
\]
where $\varphi$ is the right Frobenius-Perron eigenvector of $A$ and $\lambda_A$ the corresponding Frobenius-Perron eigenvalue, i.e., $A\varphi=\lambda_A\varphi$ and $\lambda_A>0$ and maximal. The stationary distribution $\pi_i=\hat\varphi_i\varphi_i$ where $\hat\varphi$ is the left Frobenius-Perron eigenvector of $A$, normalized so that $\sum_i\pi_i=1$. It is easy to verify that the entropy rate of the corresponding discrete-time Markov chain is in fact $\log\lambda_A$,
and that the probability of transitioning between states $i$ and $j$ in $N$ steps following any of the several possible alternative $N$-step paths is $\hat\varphi_i\varphi_j/\lambda_A^{N}$. This is a rather remarkable property, that the probability is independent of which of the alternative paths was taken! The chances that the random walker takes any of those is the same, and the entropy rate is in fact the topological entropy of the graph.
Thus, the purpose of this note is to see that there is a completely analogous continuous-time random walk with exactly the same properties.

\section{Continuous-time Ruelle Bowen random walk}

We will now resume with the construction of a continuous-time analogue of the Ruelle Bowen random walk, which is the main contribution of this work. As in the discrete-time case,
the process maximizes the entropy rate over all such continuous-time random walks. As noted, interestingly, it also equalizes the probability of traversing any path, for any specified number of transitions in its flight between any two specific nodes, in complete analogy with the standard RB random walk \cite{DelLib11}.

Consider the time-homogenous Markov jump process $X_t$ as before, with generator $Q$ and invariant measure $\pi$, and the right-continuous discretized-in-time process
\[
X_t^\Delta := X_{k\Delta}, \mbox{ for } k\Delta\leq t <(k+1)\Delta.
\]
The entropy rate of $X_t^\Delta$ is the same as that of the discrete-time Markov chain with transition probability matrix $\Pi=\exp(Q\Delta)$.  For $\Delta$ small, the entropy rate of $X_{k\Delta}$, and hence, of $X_t^\Delta$ as they carry the same amount of information, is
\[
H(X_{k\Delta}) \simeq   
-\Delta(1-\log\Delta) \sum_i\pi_i q_{ii} -\Delta \sum_{j\neq i}\pi_i q_{ij}\log q_{ij}.
\]
Clearly, the differential entropy rate $H(X_{k\Delta})/\Delta$ grows unbounded
as the time-scale resolution $\Delta$ goes to zero. A meaningful way to define a suitable rate of information for a jump process is separate the information content of the transition between nodes with that of the actual timing. Thus, we may specify at the outset the value of the average ``retention rate'' $-\sum_i\pi_iq_{ii}$ at nodes. Equivalently, we may define intstead as {\em differential entropy rate} for the jump process $X_t$
    \begin{equation}\label{eq:entropyc}
        h_\eta(Q):=-\eta\sum_{i} \pi_i q_{ii}-\sum_{i\neq j} \pi_i q_{ij}\log q_{ij},
    \end{equation} 
for any fixed choice of a value $\eta>0$. The parameter $\eta$  effectively serves as a Lagrange multiplier in the optimization that follows, that dictates a corresponding value for $-\sum_i\pi_iq_{ii}$. Whatever the interpretation, our problem is now to seek the maximum of $h_\eta(Q)$ over choices of the infinitesimal generator matrix $Q$ and in accordance with the geometry of the given graph $\cG$.

Since the graph is strongly connected, a sufficiently high power of the adjacency matrix $A$ has all of its entries strictly positive. By Frobenius-Perron theory, $A$ has a unique eigenvalue $\lambda_A$ having maximal modulus which is also real and positive. The corresponding right and left eigenvectors $\varphi$, $\hat\varphi$ can be taken to have positive entries, and also normalized so that
	\[
		\langle \varphi, \hat\varphi \rangle=\sum_{i=1}^n \varphi_i \hat\varphi_i =1.
	\]
As before, 
	\begin{equation}
		\pi^{\RB}:=\varphi\hat\varphi,
	\end{equation}
but now we also define 
	\begin{equation}\label{eq:RBc}
		Q^\RB=\diag(\varphi)^{-1} A \diag(\varphi)-\lambda_A I.
	\end{equation}
Clearly, $Q^\RB$ is a generator having $\pi^\RB$ as invariant measure, since it has non-negative off-diagonal entries, positive diagonal entries,
	\begin{align*}
		Q^\RB \one&=\diag(\varphi)^{-1} A \varphi-\lambda_A \one
		\\&=\diag(\varphi)^{-1} \lambda_A \varphi-\lambda_A \one=\zz,
	\end{align*}
and 
	\begin{align*}
		(Q^\RB)^\prime\pi^\RB&=\diag(\varphi)A' \hat\varphi-\lambda_A \pi^\RB\\ &=
		\lambda_A\diag(\varphi)\hat\varphi-\lambda_A\pi^\RB=\zz,
	\end{align*}
	where ``$\,^\prime$'' denotes transposition and $\pi^\RB$ is thought of as a column vector.
	
\begin{theorem}
The differential entropy rate $h_1(Q)$ has a unique maximum at $Q^\RB$ given above, over all infinitesimal generators that are consistent with the topology of a graph with adjacency matrix $A$. Moreover,
\[h_1(Q^\RB)=\lambda_A.
\]
\end{theorem}
{\bf Proof:} Define $r_{ij}=\pi_i Q_{ij}$ (i.e.,  for $j\neq i$ these represent ``flow rates'' from node $i$ to $j$), then maximizing $h(Q)$ reduces to maximizing 
	\[
		f(r,\pi)=-\sum_{i}r_{ii}-\sum_{i\neq j} r_{ij} a_{ij}\log \frac{r_{ij}}{\pi_i a_{ij}}
	\]
over $r, \pi$ subject to the constraints
	\[
		\sum_{j}r_{ij}=0, ~~~\mbox{and}~~~\sum_{i}r_{ij}=0.
	\]
	To see this note that $R=[r_{ij}]_{ij=1}^n$ equals $\diag(\pi)\times Q$, and the above constraints are precisely the requirements that $Q\one=\zz$ and $Q^\prime \pi =\zz$.
	
Introducing Lagrangian multipliers $\alpha, \beta$ for these two constraints gives
	\[
	\begin{split}
		\cL(r,\pi,\alpha,\beta)=-\sum_{i}r_{ii}-\sum_{i\neq j} r_{ij} a_{ij}\log \frac{r_{ij}}{\pi_i a_{ij}}\\
		+\sum_{i}\alpha_i\sum_{j}r_{ij}+\sum_{j}\beta_j\sum_{i}r_{ij}.
		\end{split}
	\]
By standard duality theory we have that the maximum of $h_1(Q)$ is bounded above by
	\[
		\min_{\alpha,\beta}\max_{r,\pi} \cL(r,\pi,\alpha,\beta).
	\]
Now choose a specific pair of $\alpha^*$ and $\beta^*$ as
	\[
		\beta_j^*=\log \varphi_j, ~~~\mbox{and}~~~ \alpha_i^*=1-\beta_i^*,
	\]	
then
	\[
		\min_{\alpha,\beta}\max_{r,\pi} \cL(r,\pi,\alpha,\beta) \le
		\max_{r,\pi} \cL(r,\pi,\alpha^*,\beta^*).
	\]
For any fixed $\pi$, we have that
	\begin{eqnarray*}
	\cL(r,\pi,\alpha^*,\beta^*) &=& \sum_{i\neq j} (-r_{ij}a_{ij}\log\frac{r_{ij}}{\pi_i a_{ij}}
	+(\alpha_i^*+\beta_j^*)r_{ij})
	\end{eqnarray*}
is a concave function of $r$. Considering first order optimality conditions we conclude that the maximizer is
	\[
		r_{ij}^*=
		\begin{cases}
		\pi_i e^{\alpha_i^*+\beta_j^*-1} & \mbox{if}~~ a_{ij}=1,\\
		0 & \mbox{if}~~ a_{ij}=0.
		\end{cases}
	\]
corresponding to the maximal value	
	\begin{eqnarray*}
	&&\sum_{i\neq j} (-r_{ij}^*a_{ij}\log\frac{r_{ij}^*}{\pi_i a_{ij}}
	+(\alpha_i^*+\beta_j^*)r_{ij}^*)
	\\&=& \sum_{i\neq j}-\pi_i e^{\alpha_i^*+\beta_j^*-1}a_{ij} \log \frac{\pi_i e^{\alpha_i^*+\beta_j^*-1}}{\pi_i a_{ij} e^{\alpha_i^*+\beta_j^*}}
	\\&=& \sum_{i\neq j}\pi_i e^{\alpha_i^*+\beta_j^*-1}a_{ij}
	\\&=& \sum_{i\neq j} \pi_i \frac{\varphi_j}{\varphi_i}a_{ij} = \lambda_A.
	\end{eqnarray*}

Therefore $h_1(Q)$ is bounded above by $\lambda_A$. On the other hand, with $Q=Q^\RB$,
	\begin{align*}
		&h_1(Q^\RB) = -\sum_{i} \pi_i q^\RB_{ii}-\sum_{i\neq j} \pi_i q^\RB_{ij}\log q^\RB_{ij}
		\\&= \lambda_A -\sum_{i\neq j} \pi_i \frac{\varphi_j}{\varphi_i}a_{ij}(\log \varphi_j
		-\log\varphi_i)
		\\&= \lambda_A-\lambda_A\sum_j \hat\varphi_j\varphi_j \log\varphi_j+
		\lambda_A\sum_i \pi_i \log\varphi_i
		\\&= \lambda_A,
	\end{align*}
and achieves the upper bound.
This completes the proof. $\Box$

The argument carries through for other values of $\eta$, in which case the solution is a scaled multiple of \eqref{eq:RBc}, corresponding to a different value for the average rate $-\sum_{i} \pi_i q_{ii}$ for staying at the current node.
\begin{cor}
The differential entropy rate $h_\eta(Q)$ has a unique maximum at $e^{\eta-1}Q^\RB$ with $Q^\RB$ given in \eqref{eq:RBc}, over all infinitesimal generators that are consistent with the topology of a graph with adjacency matrix $A$. Moreover,
\[h_\eta(e^{\eta-1}Q^\RB)=e^{\eta-1}\lambda_A.
\]
\end{cor}

The continuous-time time-homogeneous Markov jump process with infinitesimal generator $Q^\RB$ in \eqref{eq:RBc} we refer to as the continuous-time Ruelle-Bowen random walk. 
We now complete the final claim in this note that the continuous-time RB random walk shares the property, of its discrete-time counterpart, of equalizing probability across alternative paths. The proof is almost verbatim with only subtlety the need to condition on the number of steps taken in a flight between two nodes, since we are now dealing with a continuous-time jump process. But this can be easily done using the realization of the process as outlined in Section \ref{sec:realization}.

\begin{theorem} Let $X_t^\RB$ be a continuous-time RB walk on a graph $\cG$. The probability of transitioning from node $i$ to node $j$, over any specified window of time, say $[0, t_f]$,
via a sequence of precisely $N$ nodes, is independent of the particular sequence of nodes.
\end{theorem}
{\bf Proof:} 
Consider the corresponding bivariate process $(Z_k^\RB,T_k^\RB)$. The key is to observe that the 
random time $T_k$ are identically distributed with the same exponent, $-\lambda_A$, independent of $Z_k$. In fact, $Z_k$ and $T_k$ are independent. Thus,
\begin{align*}
&{\mathbb P}(T_0 = t_1, \sum_0^1T_k=t_2,\ldots, \sum_0^{N-1}T_k=t_{N},\\&\phantom{{\mathbb P}(} t_{N}\le t_f<t_{N+1},\\
&\phantom{{\mathbb P}(}X_0=i, \ldots, X_{t_f}=j)
\end{align*}
is equal to the product
\[
{\mathbb P}(t_{N}\le t_f<t_{N+1}) \times {\mathbb P}(Z_0=i,\ldots,Z_{N}=j).
\]
The former is independent of the particular path and the latter, using the facts detailed in Section \ref{sec:discrete}, is
\[
\frac{\hat\varphi_i\varphi_j}{\lambda_A^{N}}.
\]
Thus, the discrete-time analysis applies verbatim and the probability is independent of the path and only a funtion of beginning and ending nodes, and of the number of steps taken. $\Box$

The relation between the infinitesimal generator $Q^\RB$ of the continuous-time RB walk and the transition probability matrix $P$ of the discrete-time RB walk, on the same graph, is
\[
Q^\RB = \lambda_A (P^\RB - I).
\]
%

\end{document}